\documentclass{article}
\usepackage{amsthm,amsmath,amssymb}

\newcommand{\nc}{\newcommand}

\newtheorem{thm}{Theorem}[section]

\newtheorem{prop}[thm]{Proposition}
\newtheorem{lemma}[thm]{Lemma}

\newtheorem{corollary}[thm]{Corollary}

\newtheorem{definition}[thm]{Definition}

\newenvironment{defin}{\begin{definition} \rm}{\end{definition}}
\newenvironment{cor}{\begin{corollary} \rm}{\end{corollary}}

\nc{\Ext}{\operatorname{Ext}}
\nc{\NS}{\operatorname{NS}}
\nc{\Amp}{\operatorname{Amp}}
\nc{\Pic}{\operatorname{Pic}}
\nc{\Kom}{\operatorname{Kom}}
\nc{\Gr}{\operatorname{Gr}}
\nc{\Rep}{\operatorname{Rep}}
\nc{\Hom}{\operatorname{Hom}}
\nc{\Sym}{\operatorname{Sym}}
\nc{\RHom}{R\operatorname{Hom}}
\nc{\cRHom}{\operatorname{\mathcal{R}\mathcal{H}om}}
\nc{\cHom}{\operatorname{\mathcal{H}om}}
\nc{\End}{\operatorname{End}}
\nc{\Coh}{\operatorname{Coh}}
\nc{\Aut}{\operatorname{Aut}}
\nc{\Coker}{\operatornamoe{Coker}}
\nc{\coker}{\operatorname{coker}}
\nc{\Ker}{\operatorname{Ker}}
\nc{\img}{\operatorname{Im}}
\nc{\D}{\operatorname{D}}
\nc{\ch}{\operatorname{ch}}
\nc{\Stab}{\operatorname{Stab}}
\nc{\SL}{\operatorname{SL}}
\nc{\rk}{\operatorname{rk}}
\nc{\GL}{\operatorname{GL}}
\nc{\Log}{\mathop{\mathrm{Log}}}
\nc{\abs}[1]{\lvert#1\rvert}
\nc{\Cone}{\operatorname{Cone}}
\nc{\id}{\operatorname{id}}
\nc{\Li}{\operatorname{Li}}
\newcommand{\Exp}[1]{e^{#1}}
\newcommand{\I}{i}
\newcommand{\Z}{\mathbb Z}
\newcommand{\C}{\mathbb C}
\newcommand{\Dd}{\mathrm D}
\newcommand{\wt}{\widetilde}
\newcommand{\ol}{\overline}
\newcommand{\sign}{\operatorname{sign}}

\nc{\cA}{{\mathcal A}}
\nc{\cB}{{\mathcal B}}
\nc{\cC}{{\mathcal C}}
\nc{\cD}{{\mathcal D}}
\nc{\cE}{{\mathcal E}}
\nc{\cF}{{\mathcal F}}
\nc{\cG}{{\mathcal G}}
\nc{\cH}{{\mathcal H}}
\nc{\cI}{{\mathcal I}}
\nc{\cJ}{{\mathcal J}}
\nc{\cK}{{\mathcal K}}
\nc{\cL}{{\mathcal L}}
\nc{\cM}{{\mathcal M}}
\nc{\cN}{{\mathcal N}}
\nc{\cO}{{\mathcal O}}
\nc{\cP}{{\mathcal P}}
\nc{\cQ}{{\mathcal Q}}
\nc{\cR}{{\mathcal R}}
\nc{\cS}{{\mathcal S}}
\nc{\cT}{{\mathcal T}}
\nc{\cU}{{\mathcal U}}
\nc{\cV}{{\mathcal V}}
\nc{\cW}{{\mathcal W}}
\nc{\cX}{{\mathcal X}}
\nc{\cY}{{\mathcal Y}}
\nc{\cZ}{{\mathcal Z}}

\nc{\bA}{{\mathbb A}}
\nc{\bB}{{\mathbb B}}
\nc{\bC}{{\mathbb C}}
\nc{\bD}{{\mathbb D}}
\nc{\bE}{{\mathbb E}}
\nc{\bF}{{\mathbb F}}
\nc{\bG}{{\mathbb G}}
\nc{\bH}{{\mathbb H}}
\nc{\bI}{{\mathbb I}}
\nc{\bJ}{{\mathbb J}}
\nc{\bK}{{\mathbb K}}
\nc{\bL}{{\mathbb L}}
\nc{\bM}{{\mathbb M}}
\nc{\bN}{{\mathbb N}}
\nc{\bO}{{\mathbb O}}
\nc{\bP}{{\mathbb P}}
\nc{\bQ}{{\mathbb Q}}
\nc{\bR}{{\mathbb R}}
\nc{\bS}{{\mathbb S}}
\nc{\bT}{{\mathbb T}}
\nc{\bU}{{\mathbb U}}
\nc{\bV}{{\mathbb V}}
\nc{\bW}{{\mathbb W}}
\nc{\bX}{{\mathbb X}}
\nc{\bY}{{\mathbb Y}}
\nc{\bZ}{{\mathbb Z}}

\begin{document}
\title{Joyce invariants for K3 surfaces and mock theta functions}
\author{Anton Mellit and So Okada\footnote{Emails: mellit@ihes.fr and
 okada@ihes.fr, Addresses: IH\'{E}S, Le Bois-Marie, 35, route de
 Chartres, F-91440, Bures-sur-Yvette, France.}}
\maketitle

\begin{abstract}
 We will discuss Joyce invariants of stability conditions for K3
 surfaces and mock theta functions.
\end{abstract}

\section{Introduction}
Bridgeland introduced the notion of stability conditions on triangulated
categories \cite{BRD_ST}, this notion extends standard stabilities such
as Gieseker stabilities on the abelian category of coherent sheaves of 
a variety $X$, denoted by $\Coh X$, to the bounded derived category of
$\Coh X$, denoted by $\Dd(X)$.

One way to think of the notion is that it is a tool to make interesting
invariants of moduli stacks, as we have seen in the foundational work
\cite{HaNa}, in which the notion of {\it Harder-Narasimhan filtrations},
in today's term, was given birth to discuss {\it Tamagawa numbers} that
are certain volumes of moduli spaces on curves.

We would like to recall that for D-branes in superstring theory,
Douglas's work \cite{Do} on $\Pi$-stabilities motivated the notion of
stability conditions. In a Calabi-Yau variety $X$, our strings form
Riemann surfaces whose boundaries restrict to subvarieties called {\it
B-branes}, that are kinds of D-branes.  With Kontsevich's framework
\cite{Ko}, in $\Dd(X)$, the notion of $\Pi$-stabilities discusses
configuration of B-branes and its deformation, which is locally
parameterized by {\it central charges} of B-branes.  In this term, we
are taking invariants out of B-branes whose central charges align in 
the complex plane.

Now, we begin to be more specific for our paper, leaving formality a
bit out for later sections.  For stability conditions of triangulated
categories, Joyce started to extend Donaldson-Thomas invariants so that
wall-crossings of stability conditions give differential equations over
his invariants, which we call Joyce invariants.

A commutative $\bQ$ algebra $\Lambda$ containing $l$ and a {\it motivic
invariant} $I$ from the category of Artin stacks of finite type to
$\Lambda$ satisfy the following: for $I(\bC)=l$, we have
$I(\GL(n,\bC))=l^{n^2}(1-l^{-1})\cdots (1-l^{-n})$ invertible in
$\Lambda$, for quasiprojective varieties $X$ and $Y$, we have
$I(X\times Y)=I(X)I(Y)$, for a closed quasiprojective variety $Y$ in
$X$, we have $I(X)=I(Y)+I(X/Y)$, and for a quotient stack $[X/G]$ with
a special algebraic group $G$, which is a group embedded in some
$\GL(n, \bC)$ with $\GL(n,\bC)\to \GL(n, \bC)/G$ having locally trivial
fibers, we have $I([X/G])=I(X)/I(G)$.
 
 For example, some motivic invariant extends the ring structures of
 Poincar\'e or Hodge polynomials on the category of smooth projective
 varieties to our category. Generally, each motivic invariant factors
 through the ring of isomorphism classes of above quotient stacks
 $[X/G]$ \cite{Jo07b}.

 From here, we will assume that $X$ denotes an algebraic K3 surface $X$,
 and, in the stability manifold of stability conditions on $\Dd(X)$,
 $\Stab^{*}(X)$ denotes the connected component constructed by
 Bridgeland \cite{BRD_K3}.  For a stability condition $\sigma$ of
 Gieseker on $\Coh X$ or of Bridgeland on $\Dd(X)$ and {\it Mukai
 vectors} $\alpha$ in the {\it Mukai lattice} of $X$ \cite{Mu}, which is
 a nondegenerate even integer lattice, let $M^{\alpha}(\sigma)$ be
 moduli stacks of semistable objects with respect to $\sigma$; now,
 Joyce invariants $J^{\alpha}(\sigma)$ are defined with these moduli
 stacks and motivic invariants.
 
 In \cite{Jo08}, on $\Coh X$, Joyce proved that his invariants exist
 independently of the choice of Gieseker stability conditions; then, on
 $\Dd(X)$, he discussed his invariants, supposing that his invariants
 exist independently of the choice of stability conditions in
 $\Stab^{*}(X)$, which was proved by Toda \cite{To}.

 Now, with the notion of {\it numerically faithfulness} ({\it faithful}
 for short), which was introduced by the second author \cite{Ok07b}, for
 each moduli stack, the independence of the choice of stability
 conditions in $\Stab^{*}(X)$ for Joyce invariants of $\Dd(X)$ manifests
 itself as follows.

 \begin{thm}\label{prop:motivic}
  For each K3 surface $X$, Mukai vector $\alpha$ of $X$, faithful
  stability conditions $\sigma,\sigma'\in \Stab^{*}(X)$, and motivic
  invariant $I$, we have $I(M^{\alpha}(\sigma))=I(M^{\alpha}(\sigma'))$.
 \end{thm}

 Here, by \cite{Ok07b}, faithful stability conditions exist as a dense
 subset in $\Stab^{*}(X)$. So, in $\Stab^{*}(X)$, for a set of
 semistable objects with a bounded mass, by wall structures examined in
 \cite{BRD_K3}, for each Mukai vector $\alpha$ and polarization of $X$,
 we have some faithful stability condition $\sigma\in \Stab^{*}(X)$ such
 that $M^{\alpha}(\sigma)$ consists of Gieseker semistable coherent
 sheaves.

 One can check that Theorem \ref{prop:motivic} holds on any other known
 stability manifolds for Calabi-Yau surfaces such as abelian surfaces
 and minimal resolutions of surface singularities (for references of
 these stability manifolds, one can consult with the Bridgeland's survey
 \cite{BRD_06}).  Also, for some moduli stacks of stable objects and
 moduli stacks of $\mu$-semistable coherent sheaves on $X$, one can
 compare Theorem \ref{prop:motivic} with a sequence of flops in
 \cite{ArBeLi} and dimension counting in \cite{Yo}.
 
 To make explicit computation of motivic invariants, we first want to
 know our moduli stacks as moduli spaces in some details, and then
 compute isomorphism groups of objects of moduli stacks.  For example,
 for primitive Mukai vectors of positive ranks, by \cite{Yo}, moduli
 spaces of Gieseker stable coherent sheaves are deformation equivalent
 to Hilbert schemes, and they have the trivial $\bC^{*}$ isomorphism
 group for each point in the moduli stacks.

 Going beyond above primitive cases is a challenge; reasons include that
 in general moduli spaces can be singular and computing isomorphism
 groups of objects is demanding.  To explain what happens in a
 situation, for objects $E,F$ and their Mukai vectors $[E],[F]$, let
 $[E].[F]=\sum_{i}(-1)^{i}\dim \Ext^{i}(E,F)$ be the Mukai paring of
 $[E]$ and $[F]$. For objects in the moduli stack of a Mukai vector with
 non-positive self Mukai paring, in the moduli stack of a multiple of
 the Mukai vector, their direct products have a nontrivial fiber with
 some isomorphism group for each point in the fiber.
 
 On the other hand, for a Mukai vector with positive self Mukai paring,
 by \cite{Ok07b}, Theorem \ref{prop:motivic} boils down to Corollary
 \ref{cor:sph}. Let us recall Mukai vectors $\alpha$ are called {\it
 spherical}, if their self-intersections are two; in other words,
 $\alpha$ correspond to {\it spherical objects} which not only give rise
 to autoequivalences of $\Dd(X)$ \cite{ST}, but also include structure
 sheaves supported over rational curves on $X$, the structure sheaf of
 $X$, and their twists by line bundles.  Notice that each Mukai vector
 $v$ with $v.v>0$ is a multiple of a spherical Mukai vector $\alpha$.

 \begin{cor}\label{cor:sph}
  For each spherical class $\alpha$, faithful $\sigma\in \Stab^{*}(X)$,
  positive integer $n$, and motivic invariant $I$, we have
  $I(M^{n\alpha}(\sigma))=I([1/\GL(n,\bC)])=\frac{1}{l^{n^2}(1-l^{-1})\cdots
  (1-l^{-n})}$.
 \end{cor}
 
 In other words, for faithful stability conditions, we always have a
 stable spherical object for each spherical class.  As pointed out to
 the authors by Bridgeland, the existence of a stable spherical object
 of each spherical class in Corollary \ref{cor:sph} in particular gives
 another way to prove that $\Stab^{*}(X)$ is locally a bundle over the
 {\it period domain} of $X$, which consists of complexified K\"ahler
 classes of $X$ without ones that are orthogonal to spherical classes.
 
 Once we know our moduli stacks in these details, then we are able to
 compute various invariants. Indeed, after the second author discussed
 some part of the content of this paper such as Corollary
 \ref{cor:Joyce} (in the original form of Joyce invariants for some
 $\alpha$) at \cite{Ok07a} and whilst the authors were preparing this
 paper, they got notified that for standard stabilities of coherent
 sheaves of rational elliptic surfaces, Yoshioka--Nakajima computed
 their invariants \cite{YoNa}. Also, for stability conditions of
 Calabi-Yau categories of dimension three (a.k.a. 3-Calabi-Yau
 categories), Kontsevich--Soibelman discussed their invariants
 \cite{KoSo}.
 
 Here we will stick to Joyce invariants for K3 surfaces, but let us make
 some comments for our readers. Unlike invariants defined by
 Nakajima--Yoshioka, Joyce invariants involve not only arbitrary motivic
 invariants, but also correction terms of powers of $q$ based on Lie
 algebras associated to each stability condition.
 
 The invariants discussed by Kontsevich--Soibelman are (presumably)
 compatible with Joyce invariants, and they put primary emphasis on
 nontrivial wall-crossing formulas of their invariants for Calabi-Yau
 categories of dimension three.

 Now, let us go back to our case; for the Joyce invariants in Corollary
 \ref{cor:sph}, we compute as below. For the convenience of our
 formulas, we will use $q=l^{-1}$, switching between {\it Tate motive}
 and {\it Lefschetz motive}.

 \begin{cor}\label{cor:Joyce}
  $J^{n\alpha}(\sigma)=\frac{q^{n^2}}{n(1-q^n)}$.
 \end{cor}

 Here we would like to mention that we are slightly modifying the
 original formulation of Joyce invariants for K3 surfaces, as suggested
 to the authors by Zagier. Namely, in order to obtain more natural
 expressions, we omit the factor $(q^{-1}-1)$ (this is $(l-1)$ in
 \cite{Jo08}). Recall that the factor $(q^{-1}-1)$ was involved so that
 we are able to get numbers on moduli stacks of stable objects by
 replacing $q$ by one.  Instead, we take residues at $q=1$ to extend the
 notion of Euler characteristics to moduli stacks, which are not
 necessarily only of stable objects.
 
 We may regard Joyce invariants as volumes for each Mukai vector by the
 following reason.  By Theorem \ref{prop:motivic}, for each Mukai vector
 and generic choices of stability conditions, motivic invariants ignore
 the difference of moduli stacks, but unlike Joyce invariants, on
 deformations of stability conditions on stability manifolds, we do not
 know whether motivic invariants deform on moduli stacks.

 So now, we would like to take the following generating functions of
 Joyce invariants:
 \[
 J_{k} =\sum_{n>0}\frac{J^{n\alpha}(\sigma)}{n^{k}}
 = \sum_{n>0}\frac{q^{n^{2}}}{n^{k+1}(1-q^{n})}.
 \] 
 Let us point out that taking residues termwise at $q=1$ gives
 $-\zeta(k+2)$.  
 
 The generating function $J_{k}$ actually appears in the following sum
 suggested by Joyce \cite{Jo07b}.  Namely, on a stability manifold, we
 can consider the form $\sum_{\alpha\neq
 0}\frac{J^{n\alpha}(\sigma)}{Z(n\alpha)^{k}}$, which is invariant under
 autoequivalences. Also let us note that we can take the smaller form
 $\sum_{\alpha.\alpha=2}\frac{J^{n\alpha}(\sigma)}{Z(n\alpha)^{k}}$,
 which is again invariant under autoequivalences. Here we would like to
 study its building piece $J_{k}$. It is clear that cases of $k$ being
 odd give degenerated forms; so we will concentrate on cases when $k$ is
 even.

 Let us also mention that by the work of Bridgeland--Toledano-Laredo
 \cite{BrTo} and Kontsevich--Soibelman \cite{KoSo}, it has became clear
 that invariants of the moduli stacks whose images of central charges
 align make a building block to study Lie algebras associated to
 stability conditions.

 As we have seen, generating functions coming out of physics have been
 discussed with modular forms.  Now, $J_{k}$ are already some quantum
 polylogarhithms, as they are $q$-deformations of
 polylogarhithms. However, the presence of $q^{n^{2}}$ in the numerator
 does not make in particular $J_{0}$ the well-known quantum dilogarithm
 (for example, see \cite{Za}), but instead $J_{k}$ look similar to some
 of {\it mock theta functions}, which were introduced by Ramanujan
 \cite{Ra00} \cite{Ra88} and carry transformation laws similar to ones
 of theta functions.  We will pursue this view point.
 
 Let us take a quick review at mock theta functions.  The explicit
 definition on these functions was not given by Ramanujan, and this
 issue had remained for a long time.  However, quite recently, Zwegers
 in his thesis \cite{Zw} provided a way to add correction terms to the
 Ramanujan's mock theta functions to make them into {\it harmonic weak
 Maass forms of weight $\frac12$} \cite{Za_Bu}, which is explained as
 follows.

 For $\tau$ in $q=e^{ 2 \pi \I \tau}$, let $\overline{D}$ be the
 differential operator $\frac{1}{2\pi i}\frac{d}{d\overline{\tau}}$,
 $M^{k}$ be the space of meromorphic modular forms of weight $k$ for
 $k\in \frac12 \Z$ with poles only at the cusps, and $\tau=x+\I y$.
 Then, harmonic weak Maass forms of weight $k$ are real analytic modular
 forms whose derivatives with respect to $\overline{D}$ fall into the
 space $\frac{\overline{M^{2-k}}}{y^{k}}$; here, these derivatives for
 mock theta functions are called {\it shadows} \cite{Za_Bu}.
 
  Since this understanding of mock theta functions surfaced, we have
 seen achievements such as \cite{BrOn06}, \cite{BrOn07}, and
 \cite{BrOn}. Especially, Fourier coefficients of harmonic weak Maass
 forms of weight $\frac12$ played a central role, in particular, for
 solving {\it the Andrews-Dragonette Conjecture}, that is to prove an
 exact formula of Fourier coefficients of a mock theta function.
  
  Also, for an even integer $k>2$, the first author in his thesis
  \cite{Me} studied the so-called {\it higher Green's functions of
  weight $k$}, which are directly related to the harmonic weak Maass forms of
  weight $2-k$ by the {\it Maass operators} $(y^2
  \overline{D})^{\frac{k-2}2}$ and $(D+\frac{2}{4\pi\I
  y})\cdots(D+\frac{k-2}{4\pi\I y})$ with the differential operator
  $D=\frac{1}{2\pi i}\frac{d}{d\tau}$.

  Now, going back to our $J_{k}$, with certain duality, we want to
  compensate our choices of positive integers $k$.  This can be done in
  terms of differential equations, modular forms, and certain correction
  terms to $J_{k}$.  Here, differential operator $D$ may correspond to
  infinitesimal derivatives of our volumes.  Let $E_k$ and $B_{k}$ be
  the Eisenstein series and the Bernoulli numbers.  Then, we have
 \begin{align*}
  D^{k-1}J_{k-2}& = \frac{B_k}{2k}(1-E_k) - J_{-k} + \sum_{n>0}
  \frac{q^{n^2}}{ n^{1-k}}.
 \end{align*}
 
 We would like to have a duality formula which contains only modular
 forms as follows. Let us recall that in the space of modular forms of a
 given degree, Eisenstein series make distinguished basis of the
 subspace that is orthogonal to cusp forms.  Now, we take the following.
 \begin{defin}
  \[
  \cJ_{k}=\frac{B_{-k}}{2k}-\frac12\sum_{n>0} \frac{q^{n^{2}}}{n^{k+1}}+J_{k}=
  \frac{B_{-k}}{2k}+\sum_{n\neq 0} \frac{q^{n^2}}{n^{k+1}(1-q^{n})}.
  \]
 \end{defin}
  Then, this time, for positive even integers $k$, we have 
  \begin{align*}
   D^{k-1}\cJ_{k-2}+\cJ_{-k}= -\frac{B_k}{2k}E_k.
  \end{align*}
  Now, we will take $\cJ_{k}$ as granted, and study $\cJ_{-2}$ in some
  detail.  
  
  Here, 
  \begin{align*}
   \cJ_{-2}(\tau)&=-\frac{1}{24}-\frac{1}{2}\sum_{n>0}nq^{n^2}+J_{-2},
  \end{align*}
  and $\sum_{n>0}n q^{n^{2}}$ is a {\it half-theta function}.  Let
  $\theta_1(\tau) =\sum_{n\in \Z} \Exp{\pi\I n^2 \tau}$ and
  $\theta_3(\tau)=\sum_{n\in \Z} \Exp{\pi\I (n+\frac12)^2 \tau}$ be
  half-period Jacobi theta functions (at $z=0$).  Then we have the
  following.

 \begin{thm}\label{thm:mock}
  For $\SL(2, \bZ)$, with bounded growth at the cusp, there is a unique
  real analytic modular form of weight two $\tilde{g}(\tau)$ such that
  the derivative of $\tilde{g}(\tau)$ with respect to $\overline{D}$
  is $-\frac{\theta_{1}(2\tau)\overline{\theta_{1}(2\tau)}+
  \theta_{3}(2\tau)\overline{\theta_{3}(2\tau)}}{64 \pi^2
  y^{\frac{3}{2}}}$.  Now the holomorphic part of $\tilde{g}(\tau)$
  coincides with $\cJ_{-2}(\tau)$.
 \end{thm}
 
 Let us explain the words ``holomorphic part'' in Theorem
 \ref{thm:mock}; for holomorphic functions $a(\tau)$ and $b(\tau)$,
 there is a canonical way to produce a function whose derivative with
 respect to $\overline{D}$ is a function of the form $\frac{a(\tau)
 \overline{b(\tau)}}{y^k}$. It is given by the following integral
 (whenever the integral converges):
 \[
 R\left(\frac{a(\tau) \overline{b(\tau)}}{y^k}; \tau\right) := 2\pi\I a(\tau) \overline{\int_{\I\infty}^{\tau} \frac{b(z) dz}{(-\frac{\I}2(z-\bar\tau))^k}}.
 \]
 Now, the difference $\tilde{g}(\tau)-R(\overline{D} (\tilde{g}(\tau))$
 vanishes by $\overline{D}$, and we call it {\it the holomorphic part}.
 
 The story of Ramanujan's functions is parallel to Theorem
 \ref{thm:mock}, since they can be obtained as the holomorphic parts of
 certain harmonic weak Maass forms of weight $\frac12$.  Indeed, we
 prove $\cJ_{-2}(\tau)$ is in the space of mock theta functions of
 weight $\frac32$ tensored by the space $M^{\frac{1}{2}}$.
 
 Let us explain a bit more. Here, the shadow is not in the space
 $\frac{\overline{M^{\frac{1}{2}}}}{y^{\frac{3}{2}}}$, but in the
 twisted space $\frac{M^{\frac{1}{2}}\otimes
 \overline{M^{\frac{1}{2}}}}{y^{\frac{3}{2}}}$. Also, the holomorphic
 part of $\tilde{g}(\tau)$ is not a mock theta function, but a sum of
 products of ordinary theta functions of weight $\frac12$ and mock theta
 functions of weight $\frac32$, which will be derived in this paper from
 the {\it Lerch function} in \cite{Zw}. We will then be able to identify
 the Fourier coefficients of the sum to end the proof of Theorem
 \ref{thm:mock}.
 
 Now, authors are aware that we are leaving many questions open.  For
 example, we would want some understanding of moduli stacks of cases
 other than ones considered here and ${\cJ}_{k}$ for $k\neq -2$, but it
 is our impression that they rather pose fundamental questions on
 isomorphism groups of points in moduli stacks, algebras on moduli
 stacks, and mock theta functions.

 Yet, here, we investigated our cases in some detail and thank the
 Dyson's dream \cite[Section 6]{Dy}, which at some point encouraged us
 to look for mock symmetries in this context.

 \section{Definitions}
 Let us recall fundamental notions from \cite{BRD_ST}.  In this paper,
 our triangulated category $\cT$ is assumed to be $\Dd(X)$ for some K3
 surface $X$.  Let $K(\cT)$ be the {\it Grothendieck group} of
 $\cT$;i.e., $K(\cT)$ is the abelian group generated by classes of
 objects of $\cT$ such that for objects $E,F,G$ in $\cT$, we have
 $[F]=[E]+[G]$ in $K(\cT)$ whenever we have an exact triangle $E\to F\to
 G$ in $\cT$.

\subsection{Stability conditions}
A {\it stability condition} $\sigma=(Z, \cP)$ on $\cT$ consists of a
group homomorphism $Z$ from $K(\cT)$ to the complex number $\bC$ and a
family $\cP(\phi)$ of full abelian subcategories of $\cT$ indexed by
real numbers $\phi$. Each $Z$ and $\cP$ are called a {\it central
charge} and a {\it slicing}. They need to satisfy the following
compatibilities.

\begin{itemize}
 \item If for some $\phi\in \bR$, $E$ is a nonzero object in
       $\cP(\phi)$, then for some positive real number $m(E)$, called
       {\it mass} of $E$, we have
       $Z(E)=m(E)\exp(i\pi\phi)$.
 \item For each real number $\phi$, we have $\cP(\phi+1)=\cP(\phi)[1]$.
 \item For real numbers $\phi_{1}>\phi_{2}$ and objects $A_{i}\in
       \cP(\phi_{i})$, we have $\Hom_{\cT}(A_{1}, A_{2})=0$.
 \item For any nonzero object $E \in \cT$, there exist real numbers
       $\phi_{1}> \cdots >\phi_{n}$ and objects $A_{i}\in \cP(\phi_{i})$
       such that there exists a sequence of exact triangles $E_{i-1}\to
       E_{i}\to A_{i}$ with $E_{0}=0$ and $E$.
\end{itemize}
The sequence above is called the {\it Harder-Narasimhan filtration}
({\it HN-filtration} for short) of $E$. The HN-filtration of any object
is unique up to isomorphisms. For each $\phi\in \bR$, nonzero objects in
$\cP(\phi)$ are called {\it semistable} with {\it phase} $\phi$.  If
moreover a semistable object in $\cP(\phi)$ has only the trivial
Jordan-H\"older filtration in $\cP(\phi)$, then it is called {\it
stable}.

We will assume that our central charge $Z$ factors through the map
$[E]\in K(\cT)\mapsto \ch(E)\sqrt{X}$, which is the {\it Mukai vector}
of $E$ in the {\it Mukai lattice} of $X$. For Mukai vectors $v,w$, let
$v.w$ be {\it the Mukai paring}.

A stability condition $\sigma=(Z,\cP)$ is called {\it numerically
faithful} \cite[Definition 3.1]{Ok07b} ({\it faithful} for short), if for
each real number $r$, we have a primitive Mukai vector $v$ such that for
each semistable object $E$ of the phase $r$, $[E]$ is a sum of
$v$. Here, by \cite[Proposition 8.3]{BRD_K3}, the connected component
$\Stab^{*}(X)$ satisfies the assumption of \cite[Lemma 3.1]{Ok07b}.  So
faithful stability conditions are dense in $\Stab^{*}(X)$.

For a real number $r$, let $\cP(r-1,r]$ be the extension-closed full
subcategory consisting of semistable objects whose phases are in the
interval $(r-1,r]$, and $C(r)$ be the Mukai vectors of the objects in
$\cP(r-1,r]$.  Then, for each stability condition $\sigma=(Z,\cP)\in
\Stab^{*}(X)$, Mukai vector $\alpha$, and the real number $r$ such that
$Z(\alpha)\in \bR_{>0} e^{i\pi r}$, we define the Joyce invariant
$J^{\alpha}(\sigma)$ to be $\sum_{n=1}^{\infty}\sum_{\alpha_1+\cdots +
\alpha_n=\alpha, \alpha_i\in C(r)} q^{\sum_{j>i} \alpha_j.\alpha_i}
\frac{(-1)^{n-1}}{n}\Pi_{i=1}I(M^{\alpha_i}(\sigma))$ \cite[Definition
6.22]{Jo08}, \cite[Definition 5.9]{To} (let us recall that as explained
in the introduction, we let $q=l^{-1}$ and omit $(l-1)$ from their
original definitions).

\subsection{Modular forms}

Let us recall the definition and properties of the Dedekind eta function
(we denote $q=\Exp{2\pi\I\tau})$, $\tau$ belongs to the upper half plane.
\[
\eta(\tau) 
= \Exp{\frac{\pi\I\tau}{12}} \prod_{n=1}^\infty (1-\Exp{2\pi\I n\tau})
= q^{1/24} (1-q-q^2+q^5+q^7\cdots).
\]
The eta functions transforms like a modular form of weight $\frac12$:
\[
\eta(\tau+1) = \Exp{\frac{\pi\I}{12}} \eta(\tau),\qquad
\eta\left(\frac{-1}{\tau}\right) = \Exp{\frac{\pi\I}4} \sqrt{\tau} \eta(\tau).
\]
We will need the following identity:
\begin{equation}\label{eq:my3}
\eta\left(\frac{\tau}2\right) \eta\left(\frac{1+\tau}2\right) \eta\left(2\tau\right)
= \Exp{\frac{\pi\I}{24}} \eta(\tau)^3.
\end{equation}

Next we recall the half-period Jacobi theta functions (at $z=0$), note
that we slightly changed the indexing:
\begin{align*}
\theta_1(\tau) 
&= \theta_{00}(0; \tau) = \sum_{n\in \Z} \Exp{\pi\I n^2 \tau} 
&= 1 + 2 q^\frac12 + 2 q^2 + 2 q^\frac92 + 2 q^8 + \cdots, \\
\theta_2(\tau)
&= \theta_{01}(0;\tau) = \sum_{n\in \Z} (-1)^n \Exp{\pi\I n^2 \tau} 
&= 1 - 2 q^\frac12 + 2 q^2 - 2 q^\frac92 + 2 q^8 + \cdots,\\
\theta_3(\tau)
&= \theta_{10}(0;\tau) = \sum_{n\in \Z} \Exp{\pi\I (n+\frac12)^2 \tau} 
&= 2 q^\frac18 + 2 q^\frac98 + 2 q^\frac{25}8 + 2 q^\frac{49}8 + \cdots.
\end{align*}

The theta functions can be expressed in terms of the eta function in the following way:
\begin{equation}\label{eq:my5}
\theta_1(\tau) = \Exp{-\frac{\pi\I}{12}}\frac{\eta\left(\frac{1+\tau}2\right)^2} {\eta(\tau)},\qquad
\theta_2(\tau) = \frac{\eta\left(\frac{\tau}2\right)^2} {\eta(\tau)},\qquad
\theta_3(\tau) = 2 \frac{\eta\left(2\tau\right)^2} {\eta(\tau)}.
\end{equation}
We know their transformation properties:
\begin{align*}
\theta_1(\tau+1) &= \theta_2(\tau), &\qquad 
\theta_1\left(-\frac{1}{\tau}\right) &= \Exp{-\frac{\pi\I}4} \sqrt{\tau} \theta_1(\tau), \\
\theta_2(\tau+1) &= \theta_1(\tau), &\qquad
\theta_2\left(-\frac{1}{\tau}\right) &= \Exp{-\frac{\pi\I}4} \sqrt{\tau} \theta_3(\tau), \\
\theta_3(\tau+1) &= \Exp{\frac{\pi\I}4} \theta_3(\tau), &\qquad
\theta_3\left(-\frac{1}{\tau}\right) &= \Exp{-\frac{\pi\I}4} \sqrt{\tau} \theta_2(\tau).
\end{align*}
In particular, they are modular forms for the group $\Gamma(2)$.

We also need the classical Eisenstein series of weight $2$ for $\SL(2,\Z)$:
\[
E_2(\tau) = 24\frac{\eta'(\tau)}{\eta(\tau)} = 1 - 24 \sum_{k,n=1}^{\infty} k q^{nk}
= 1 - 24 q - 72 q^2 - 96 q^3 - \cdots.
\]

The function $E_2$ is not a modular form, but is quasi-modular form. The fourth powers of the theta functions are Eisenstein series for $\Gamma(2)$ and we have the following relations:
\begin{align}\label{eq:my8}
E_2\left(\frac{1+\tau}2\right) - 2 E_2(\tau) &= \theta_1^4(\tau) - 2 \theta_2(\tau)^4,\\
E_2\left(\frac{\tau}2\right) - 2 E_2(\tau) &= -2 \theta_1^4(\tau) + \theta_2(\tau)^4,\\
4 E_2\left(2 \tau \right) - 2 E_2(\tau) &= \theta_1^4(\tau) + \theta_2(\tau)^4,\\
\theta_1^4(\tau) &= \theta_2^4(\tau) + \theta_3^4(\tau).
\end{align}

\subsection{The Lerch function}

Having introduced some classical modular forms, we turn to the thesis of
Zwegers \cite{Zw}. In this thesis we find the following definition of
the Lerch function:
\[
\mu(u, v; \tau) = \frac{\Exp{\pi\I u}} {\theta(v;\tau)} \sum_{n\in\Z} 
\frac{(-1)^n \Exp{\pi\I(n^2+n)\tau+2\pi\I n v}} {1 - \Exp{2\pi\I n \tau + 2\pi\I u}} \qquad (u,v\in \C\setminus (\Z\tau+\Z)).
\]
The definition of the theta function he uses is the following one:
\[
\theta(z;\tau) = \sum_{\nu\in \frac12+\Z} \Exp{\pi\I \nu^2\tau + 2\pi\I\nu(z+\frac12)}.
\]
Note the following symmetry:
\[
\theta(z+1;\tau) = \theta(-z;\tau) = -\theta(z;\tau), \qquad \theta(z+\tau;\tau) = -\Exp{-\pi\I\tau-2\pi\I z} \theta(z;\tau).
\]
The theta functions $\theta_1$, $\theta_2$ and $\theta_3$ are related to $\theta$ in the following way:
\begin{align*}
\theta\left(\frac{\tau}2;\tau\right) &= 
    -\I \Exp{-\frac{\pi\I\tau}4} \theta_2(\tau),\\
\theta\left(\frac{1+\tau}2;\tau\right) &=
    -\Exp{-\frac{\pi\I\tau}4} \theta_1(\tau),\\
\theta\left(\frac12;\tau\right) &=
    - \theta_3(\tau).
\end{align*}
Moreover we have
\[
\theta(0;\tau)=0,\qquad \frac{d}{2\pi\I d s}\bigg\vert_{s=0} \theta(s;\tau) = \I \eta^3(\tau).
\]

Zwegers found a way to add a correction term to $\mu$ so that the new
function $\wt\mu$ has good transformation properties. Namely, he defines
\[
\wt\mu(u,v;\tau) = \mu(u,v;\tau) + \frac{\I}2 R(u-v;\tau),
\]
where
\[
R(u;\tau) = \sum_{\nu\in\frac12+\Z} \left\{\sign(\nu) - E\left( \left(\nu+\frac{\Im u} y \right) \sqrt{2y}\right) \right\} (-1)^{\nu-\frac12} \Exp{-\pi\I\nu^2\tau-2\pi\I\nu u}.
\]
Here $y = \Im \tau$ and $E$ is the function
\[
E(z) = 2\int_0^z \Exp{-\pi t^2} dt = 1-\mathrm{erfc}(z\sqrt{\pi}).
\]
The result of Zwegers is the following transformation properties of $\wt\mu$:
\begin{thm}\cite[Theorem 1.11]{Zw}\label{thm:my1}
The function $\wt\mu$ satisfies 
\[
\wt\mu(u,v;\tau)=\wt\mu(v,u;\tau)=\wt\mu(-u,-v;\tau),
\]
and 
\begin{align*}
\wt\mu(u,v;\tau+1) &= \Exp{-\frac{\pi\I}4} \wt\mu(u,v;\tau), &\qquad
\wt\mu\left(\frac{u}{\tau}, \frac{v}{\tau}; -\frac{1}{\tau}\right) &= - \Exp{-\frac{\pi\I}4 - \frac{\pi\I (u-v)^2}{\tau}} \sqrt{\tau} \wt\mu(u,v;\tau)\\
\wt\mu(u+1,v;\tau) &= - \wt\mu(u,v;\tau), &\quad 
\wt\mu(u+\tau,v;\tau) &= -\Exp{2\pi\I(u-v)+\pi\I\tau} \wt\mu(u,v).
\end{align*}
\end{thm}

Here is a list of properties that the functions $R$ and $\mu$ satisfy separately:
\begin{prop}\cite[Propositions 1.4 and 1.9]{Zw}
The functions $\mu$ and $R$ satisfy
\[
\mu(u,v;\tau)=\mu(v,u;\tau)=\mu(-u,-v;\tau), \qquad R(-z;\tau) = R(z;\tau),
\]
and we have
\[
\mu(u+1,v;\tau) = - \mu(u,v;\tau), \qquad R(z+1;\tau)=-R(z;\tau).
\]
\end{prop}

We also mention one last property which we will use:
\begin{prop}\cite[Proposition 1.4 and Theorem 1.11]{Zw}\label{prop:my3}
Both the function $\mu$ and $\wt\mu$ (if you plug it in place of $\mu$) satisfy
\[
\mu(u+z, v+z;\tau) - \mu(u, v;\tau) = \frac{\I \eta^3(\tau) \theta(u+v+z;\tau) \theta(z;\tau)}
{\theta(u;\tau) \theta(v;\tau) \theta(u+z;\tau) \theta(v+z;\tau)}
\]
for $u,v,u+z,v+z \notin \Z+\tau\Z$.
\end{prop}

\section{Proofs}

Let us prove Theorem \ref{prop:motivic}.
\begin{proof}
 For faithful stability conditions $\sigma$, in terms of Mukai vectors,
 $J^{\alpha}(\sigma)$ admit unique expressions.  Since by \cite[Theorem
 1.5]{To}, we have $J^{\alpha}(\sigma)=J^{\alpha}(\sigma')$ for any
 $\alpha$, especially for primitive ones, the statement follows.
\end{proof}

 In terms of faithful stability conditions over integer lattices, for
 invariants of moduli stacks of aligned central charges, Theorem
 \ref{prop:motivic} is a general feature of their deformation invariance
 on stability manifolds.
 
 We will prove Corollary \ref{cor:sph}. Now, an object $E\in \cT$ is
 called {\it spherical} if $\Ext^{i}(E,E)=\bC$ for $i=0, 2$ and
 $\Ext^{i}(E,E)=0$ for else; spherical classes are Mukai vectors of
 spherical objects.

\begin{proof}
 For the case when $\alpha$ is with a nonzero rank, by \cite[Theorem
 0.1(1)]{Yo} and \cite[Proposition 14.2]{BRD_K3}, for some faithful
 $\sigma\in \Stab^{*}(X)$, we have a stable spherical object whose 
 class is $\alpha$. So, by \cite[Proposition 4.9]{Ok07b}, the 
 statement follows.  For other cases, by \cite[Lemma 25]{Fr}, the first 
 Chern class of $\alpha$ is either effective or anti-effective. So, by
 replacing $\alpha$ with $-\alpha$, if necessarily, one recalls that
 some coherent sheaf $E$ with $[E]=\alpha$ is Gieseker semistable.
 Then, by \cite[Theorem 6.6]{To}, the statement follows.
\end{proof}

Let us prove Corollary \ref{cor:Joyce}.

\begin{proof}
 Since $\alpha.\alpha=2$, by choosing $\sigma$ to be faithful, we have
 that for positive integers $k_{i}$, $J^{n\alpha}(\sigma)$ is equal to
 $\sum_{m=1}^{\infty}\sum_{k_{1}+\cdots +k_m=n} q^{\sum_{i>j} 2 k_i
 k_j}\frac{(-1)^{n-1}}{n} \Pi_{i=1}^{n} \frac{1}{I(\GL(k_{i}, \bC))}$.
 Since $\sum_{i>j} 2 k_i k_j=(\sum k_i)^2- \sum k_i^2=n^2-\sum k_i^2$,
 we have that $J^{n\alpha}(\sigma)$ is equal to
 $q^{n^2}\sum_{m=1}^{\infty} \sum_{k_{1}+\cdots
 +k_{m}=n}\frac{(-1)^{n-1}}{n}\Pi_{i=1}^{n}
 \frac{q^{-k_i^2}}{I(\GL(k_{i}, \bC))}$.
 
 Let $F(x)=\sum_{m\geq 0} \frac{q^{-m^2}}{I(\GL(m,\bC))}x^{m}$. Then we
 have $F(x)-F(q x)=x F(x)$, and $J^{n\alpha}$ is the $n$-th coefficient
 of $q^{n^2} \sum \frac{(-1)^{n-1}}{n} (F(x)-1)=q^{n^2} \log F(x)$.
 Since $\log F(x)+\log(1-x)=\log F(q x)$, the $n$-th coefficient of
 $\log F(x)$ is $\frac{1}{n (1-q^{n})}$. So the statement follows.
\end{proof}

The rest of this section is devoted to the proof of Theorem
\ref{thm:mock}.  The plan is to see the existence of a function with the
holomorphic part being $\cJ_{-2}$ and good transformation
properties. Now, the first clue is to notice that $\cJ_{-2}$ looks
similar to $\mu$, which is the holomorphic part of $\tilde{\mu}$, but to
be precise, we will here derive several functions from $\mu$ and
subsequently modify them with theta functions.

Let us study behavior of the functions $\mu, \wt\mu, R$ at the ``points
of order two''. The values at these points are not interesting since we
have
\begin{prop}
 \[
 \wt\mu\left(\frac12, \frac{\tau}2;\tau\right) = \wt\mu\left(\frac12, \frac{1+\tau}2;\tau\right) = \wt\mu\left(\frac{\tau}2, \frac{1+\tau}2;\tau\right) = 0.
 \]
\end{prop}
\begin{proof}
 We simply take the definition of $\mu$ and $R$ above and use the
 following trick. For example, in the case of $\wt\mu\left(\frac12,
 \frac{\tau}2\right)$ the trick is to write
\begin{align*}
\sum_{n\in\Z} \frac{(-1)^n \Exp{\pi\I(n^2+2n)\tau}}{1+\Exp{2\pi\I n\tau}}
&=\frac12 \left(\sum_{n\in\Z} \frac{(-1)^n \Exp{\pi\I(n^2+2n)\tau}}{1+\Exp{2\pi\I n\tau}}
+ \sum_{n\in\Z} \frac{(-1)^n \Exp{\pi\I(n^2-2n)\tau}}{1+\Exp{-2\pi\I n\tau}} \right)\\
&=\frac12 \sum_{n\in\Z} (-1)^n \Exp{\pi\I n^2\tau} = \frac{\theta_2(\tau)}2.
\end{align*}

 Therefore for $\wt\mu\left(\frac12, \frac{\tau}2\right)$ we obtain
\[
\mu\left(\frac12, \frac{\tau}2;\tau\right) = - \frac{\Exp{\frac{\pi\I\tau}4}}2.
\]
A trick similar to the one used above gives
\[
R\left(\frac{1-\tau}2;\tau\right) = -\I \Exp{\frac{\pi\I\tau}4}.
\]
Thus $\wt\mu\left(\frac12, \frac{\tau}2;\tau\right) = 0$. The other cases are similar with
\begin{align*}
\mu\left(\frac12, \frac{1+\tau}2;\tau\right) &= -\frac{\I\Exp{\frac{\pi\I\tau}4} }2, &\qquad
\mu\left(\frac{\tau}2, \frac{1+\tau}2;\tau\right) &= 0,\\
R\left(-\frac{\tau}2;\tau\right) &= \Exp{\frac{\pi\I\tau}4}, &\qquad
R\left(-\frac12;\tau\right) &=0.
\end{align*}
\end{proof}

Because of the last proposition the derivatives of $\wt\mu$ at the points of order $2$ should have nice transformation properties. Namely, we define
\begin{align*}
\mu_1(\tau) &= \frac{d}{2\pi\I d s}\bigg\vert_{s=0} \wt\mu\left(\frac12, \frac{\tau}2 + s;\tau \right), &\quad
\mu'_1(\tau) &= \frac{d}{2\pi\I d s}\bigg\vert_{s=0} \wt\mu\left(\frac12 + s, \frac{\tau}2;\tau \right),\\
\mu_2(\tau) &= \frac{d}{2\pi\I d s}\bigg\vert_{s=0} \wt\mu\left(\frac12, \frac{1+\tau}2 + s;\tau \right), &\quad
\mu'_2(\tau) &= \frac{d}{2\pi\I d s}\bigg\vert_{s=0} \wt\mu\left(\frac12 + s, \frac{1+\tau}2;\tau \right),\\
\mu_3(\tau) &= \frac{d}{2\pi\I d s}\bigg\vert_{s=0} \wt\mu\left(\frac{\tau}2, \frac{1+\tau}2 + s;\tau \right), &\quad
\mu'_3(\tau) &= \frac{d}{2\pi\I d s}\bigg\vert_{s=0} \wt\mu\left(\frac{\tau}2 + s, \frac{1+\tau}2;\tau \right).
\end{align*}

\begin{prop}
We have
\begin{align*}
\mu_1(\tau) + \mu'_1(\tau) &= - \frac{\Exp{\frac{\pi\I\tau}4} \theta_1(\tau)^3}4,\\
\mu_2(\tau) + \mu'_2(\tau) &= - \frac{\I \Exp{\frac{\pi\I\tau}4} \theta_2(\tau)^3}4,\\
\mu_3(\tau) + \mu'_3(\tau) &= -\frac{\theta_3(\tau)^3}4.
\end{align*}
\end{prop}
\begin{proof}
Using Proposition \ref{prop:my3} we obtain
\begin{align*}
\mu_1(\tau) + \mu'_1(\tau) 
&= \frac{d}{2\pi\I d s}\bigg\vert_{s=0} \left( \wt\mu\left(\frac12, \frac{\tau}2 + s \right) 
- \wt\mu\left(\frac12 - s, \frac{\tau}2 \right) \right)\\
&= \frac{d}{2\pi\I d s}\bigg\vert_{s=0} 
\frac{\I\eta(\tau)^3 \theta\left(\frac{1+\tau}2; \tau\right) \theta\left(s;\tau\right)}
{\theta\left(\frac12-s;\tau\right) \theta\left(\frac12;\tau\right) \theta\left(\frac{\tau}2-s;\tau\right)
\theta\left(\frac{\tau}2;\tau\right)}\\
&= \frac{-\eta(\tau)^6 \theta\left(\frac{1+\tau}2; \tau\right)}
{ \theta\left(\frac12;\tau\right)^2 \theta\left(\frac{\tau}2;\tau\right)^2} 
= \frac{-\eta^6(\tau) \Exp{\frac{\pi\I\tau}4} \theta_1(\tau)}{\theta_2(\tau)^2\theta_3(\tau)^2} .
\end{align*}
We have (using \eqref{eq:my5} and \eqref{eq:my3})
\[
\theta_1(\tau) \theta_2(\tau) \theta_3(\tau) = 2 \eta(\tau)^3.
\]
Therefore
\[
\mu_1(\tau) + \mu'_1(\tau)  = - \frac{\Exp{\frac{\pi\I\tau}4} \theta_1(\tau)^3}4.
\]
Similarly,
\begin{align*}
\mu_2(\tau) + \mu'_2(\tau) 
&=\frac{-\eta(\tau)^6 \theta\left(\frac{\tau}2+1; \tau\right)}
{ \theta\left(\frac12;\tau\right)^2 \theta\left(\frac{1+\tau}2;\tau\right)^2} 
= \frac{-\I \eta^6(\tau) \Exp{\frac{\pi\I\tau}4} \theta_1(\tau)}{\theta_2(\tau)^2\theta_3(\tau)^2}
=- \frac{\I \Exp{\frac{\pi\I\tau}4} \theta_2(\tau)^3}4,\\
\mu_3(\tau) + \mu'_3(\tau) 
&=\frac{-\eta(\tau)^6 \theta\left(\frac12+\tau; \tau\right)}
{ \theta\left(\frac{\tau}2;\tau\right)^2 \theta\left(\frac{1+\tau}2;\tau\right)^2} 
= -\frac{\eta(\tau)^6 \theta_3(\tau)}{\theta_1(\tau)^2 \theta_2(\tau)^2}
= -\frac{\theta_3(\tau)^3}4.
\end{align*}
\end{proof}

Now we are ready to formulate the transformation properties of $\mu_i$.
\begin{prop}
We have
\begin{align*}
\mu_1(1+\tau) &= \Exp{-\frac{\pi\I}4} \mu_2(\tau), \qquad
\mu_2(1+\tau) = -\Exp{-\frac{\pi\I}4} \mu_1(\tau),\\
\mu_3(1+\tau) &= \Exp{-\frac{\pi\I}4} \left(\frac{\theta_3(\tau)^3}4 + \mu_3(\tau)\right),\\
\mu_1\left(-\frac1\tau\right) &= 
\Exp{\frac{\pi\I}4-\pi\I\frac{1+\tau^2}{4\tau}} \tau^{\frac32} \left(\frac{\Exp{\frac{\pi\I\tau}4}\theta_1(\tau)^3} 4 + \mu_1(\tau)\right), \\
\mu_2\left(-\frac1\tau\right) &= 
\Exp{-\frac{\pi\I}4 - \frac{\pi\I}{4\tau}} \tau^{\frac32}\mu_3(\tau), \qquad
\mu_3\left(-\frac1\tau\right) =
-\Exp{-\frac{\pi\I}4 - \frac{\pi\I\tau}4} \tau^{\frac32}\mu_2(\tau).
\end{align*}
\end{prop}
\begin{proof}
The proof of the first three equations goes by applying the operator $\frac{d}{2\pi\I d s}\big\vert_{s=0}$ to both sides of the following equations obtained from Theorem \ref{thm:my1}.
\begin{align*}
\wt\mu\left(\frac12,\frac{1+\tau}2+s;\tau+1\right) 
&= \Exp{-\frac{\pi\I}4} \wt\mu\left(\frac12, \frac{1+\tau}2+s;\tau\right),\\
\wt\mu\left(\frac12,\frac{\tau}2+1+s;\tau+1\right)
&= -\Exp{-\frac{\pi\I}4} \wt\mu\left(\frac12, \frac{\tau}2+s;\tau\right),\\
\wt\mu\left(\frac{1+\tau}2, \frac{\tau}2 + 1 + s;\tau+1\right)
&= -\Exp{-\frac{\pi\I}4} \wt\mu\left(\frac{\tau}2 + s, \frac{1+\tau}2;\tau\right).
\end{align*}

Similarly, for the last three equations we use
\begin{align*}
\wt\mu\left(\frac12,-\frac{1}{2\tau}+s;-\frac{1}{\tau}\right)
&= -\Exp{-\frac{\pi\I}4 - \frac{\pi\I(\tau+1-2 s \tau)^2}{4\tau}} \sqrt{\tau} \wt\mu\left(\frac{\tau}2, -\frac12 + s\tau;\tau\right),\\
\wt\mu\left(\frac12, \frac{\tau-1}{2\tau}+s; -\frac{1}{\tau}\right)
&= -\Exp{-\frac{\pi\I}4 - \frac{\pi\I(1-2s\tau)^2}{4\tau}} \sqrt{\tau} \wt\mu\left(\frac{\tau}2, \frac{\tau-1}2 + s\tau; \tau\right),\\
\wt\mu\left(-\frac{1}{2\tau}, \frac{\tau-1}{2\tau}+s; -\frac{1}{\tau}\right)
&= -\Exp{-\frac{\pi\I}4 - \frac{\pi\I(\tau+2s\tau)^2}{4\tau}} \sqrt{\tau} \wt\mu\left(-\frac12, \frac{\tau-1}2 + s\tau; \tau\right).
\end{align*}
\end{proof}

Looking at the proposition above it is clear that we should consider the following three functions
\begin{align*}
\wt h_1(\tau) &= \Exp{-\frac{\pi\I\tau}4} \mu_1(\tau) + \frac{\theta_1(\tau)^3}8, \\
\wt h_2(\tau) &= -\I \Exp{-\frac{\pi\I\tau}4} \mu_2(\tau) + \frac{\theta_2(\tau)^3}8, \\
\wt h_3(\tau) &= -\mu_3(\tau) - \frac{\theta_3(\tau)^3}8.
\end{align*}

Then we have
\begin{align*}
\wt h_1(\tau+1) &= \wt h_2(\tau), &\qquad 
\wt h_1\left(-\frac{1}{\tau}\right) &= \Exp{\frac{\pi\I}4} \tau^{\frac32} \wt h_1(\tau), \\
\wt h_2(\tau+1) &= \wt h_1(\tau), &\qquad
\wt h_2\left(-\frac{1}{\tau}\right) &= \Exp{\frac{\pi\I}4} \tau^{\frac32} \wt h_3(\tau), \\
\wt h_3(\tau+1) &= \Exp{-\frac{\pi\I}4} \wt h_3(\tau), &\qquad
\wt h_3\left(-\frac{1}{\tau}\right) &= \Exp{\frac{\pi\I}4} \tau^{\frac32} \wt h_2(\tau).
\end{align*}

Next we need to find the Fourier expansions of $\wt h_i$. We would like
to have them similar to the decomposition $\wt\mu = \mu + \frac{\I}2
R$. Therefore we compute
\begin{prop}\label{prop:my7}
\begin{align*}
\frac{d}{2\pi\I d s}\bigg\vert_{s=0} R\left(\frac{1-\tau}2 - s; \tau \right) 
 &= -\I \Exp{\frac{\pi\I\tau}4} \left( \sum_{n\in\Z} |n| \beta(2y n^2)  \Exp{ - \pi\I n^2\tau}
+\frac12 - \frac{\ol{\theta_1(\tau)}}{\pi\sqrt{2y}}\right), \\
\frac{d}{2\pi\I d s}\bigg\vert_{s=0} R\left(-\frac{\tau}2 - s; \tau \right) 
 &= \Exp{\frac{\pi\I\tau}4} \left( \sum_{n\in\Z} (-1)^n |n| \beta(2y n^2)  \Exp{ - \pi\I n^2\tau}
+\frac12 - \frac{\ol{\theta_2(\tau)}}{\pi\sqrt{2y}}\right), \\
\frac{d}{2\pi\I d s}\bigg\vert_{s=0} R\left(-\frac12 - s; \tau \right) 
 &= \I \left( \sum_{\nu\in\Z+\frac12} |\nu| \beta(2y\nu^2) \Exp{ -\pi\I\nu^2\tau} 
 - \frac{\ol{\theta_3(\tau)}}{\pi\sqrt{2y}} \right).
\end{align*}
\end{prop}
\begin{proof}
Differentiating term by term gives
\begin{multline*}
\frac{d}{2\pi\I d s}\bigg\vert_{s=0} R (u-s;\tau) 
\\= \sum_{\nu\in\frac12+\Z} \nu \left\{\sign(\nu) - E\left(\left(\nu+\frac{\Im u} y\right) \sqrt{2y} \right) \right\} (-1)^{\nu-\frac12} \Exp{-\pi\I\nu^2\tau - 2\pi\I\nu u}\\
-\sum_{\nu\in\frac12 + \Z} \frac{1}{\pi \sqrt{2y}} (-1)^{\nu-\frac12} \Exp{-2\pi y \left(\nu+\frac{\Im u} y\right)^2 -\pi\I\nu^2\tau - 2\pi\I\nu u}.
\end{multline*}
Therefore for the first case we obtain
\begin{multline*}
-\I \Exp{\frac{\pi\I\tau}4} \sum_{n\in\Z} \left(n+\frac12\right) \left(\sign\left(n+\frac12\right) - E(n\sqrt{2y})\right) \Exp{ - \pi\I n^2\tau}\\
+ \frac{\I\Exp{\frac{\pi\I\tau}4}}{\pi \sqrt{2y}} \sum_{n\in\Z} \Exp{-\pi\I n^2 \bar \tau}.
\end{multline*}
The first summand can be transformed into
\[
-\I \Exp{\frac{\pi\I\tau}4} \sum_{n\in\Z} \left(n+\frac12\right)(\sign(n) - E(n\sqrt{2y})) \Exp{ - \pi\I n^2\tau}
-\frac{\I \Exp{\frac{\pi\I\tau}4}}2.
\]
Using the function $\beta$,
\[
\beta(x) = \sum_x^\infty t^{-\frac12}\Exp{-\pi t}dt = 1-E(\sqrt{x}) = \mathrm{erfc}\;(\sqrt{\pi x}),
\]
we obtain
\begin{multline*}
\frac{d}{2\pi\I d s}\bigg\vert_{s=0} R\left(\frac{1-\tau}2 - s; \tau \right)
\\=
-\I \Exp{\frac{\pi\I\tau}4} \sum_{n\in\Z} \left(n+\frac12\right) \sign(n) \beta(2y n^2)  \Exp{ - \pi\I n^2\tau}
-\frac{\I \Exp{\frac{\pi\I\tau}4}}2 + \frac{\I\Exp{\frac{\pi\I\tau}4}}{\pi\sqrt{2y}} \ol{\theta_1(\tau)},
\end{multline*}
and the final result easily follows from this formula.

Analogously, in the second case we obtain
\begin{multline*}
\frac{d}{2\pi\I d s}\bigg\vert_{s=0} R\left(-\frac{\tau}2 - s; \tau \right)
\\=
\Exp{\frac{\pi\I\tau}4} \sum_{n\in\Z} (-1)^n \left(n+\frac12\right) \sign(n) \beta(2y n^2)  \Exp{ - \pi\I n^2\tau}
+\frac{\Exp{\frac{\pi\I\tau}4}}2 - \frac{\Exp{\frac{\pi\I\tau}4}}{\pi\sqrt{2y}} \ol{\theta_2(\tau)},
\end{multline*}

In the third case the result follows right from the following formula:
\[
\frac{d}{2\pi\I d s}\bigg\vert_{s=0} R\left(-\frac12 - s; \tau \right)
=\I \sum_{\nu\in\Z+\frac12} |\nu| \beta(2 y \nu^2) \Exp{-\pi\I \nu^2\tau}
- \I \sum_{\nu\in\Z+\frac12} \frac{\Exp{-\pi\I \nu^2\bar\tau}}{\pi\sqrt{2y}}.
\]
\end{proof}

It remains to differentiate the function $\mu$.
\begin{prop}\label{prop:my8}
The corresponding derivatives of $\mu$ are given by
\begin{align*}
    \frac{d}{2\pi\I d s}\bigg\vert_{s=0} &\mu\left(\frac12, \frac{\tau}2+s;\tau\right)\\
& = -\frac{\Exp{\frac{\pi\I\tau}4}}{24 \theta_1(\tau)} \left(
-2 + 6 \theta_1(\tau) + 3 \theta_1^4(\tau)  - E_2(\tau)
 + 48 \sum_{n=1}^\infty \frac{\Exp{\pi\I (n^2+2n)\tau}}{(1-\Exp{2\pi\I n\tau})^2}\right),\\
    \frac{d}{2\pi\I d s}\bigg\vert_{s=0} &\mu\left(\frac12, \frac{1+\tau}2+s;\tau\right)\\
& = -\I \frac{\Exp{\frac{\pi\I\tau}4}}{24 \theta_2(\tau)} \left(
-2 + 6 \theta_2(\tau) + 3 \theta_2^4(\tau) - E_2(\tau)
+48 \sum_{n=1}^\infty \frac{(-1)^n \Exp{\pi\I (n^2+2n)\tau}}{(1-\Exp{2\pi\I n\tau})^2} \right),\\
    \frac{d}{2\pi\I d s}\bigg\vert_{s=0} &\mu\left(\frac{\tau}2, \frac{1+\tau}2+s;\tau\right)\\
& =\frac{1}{24 \theta_3(\tau)} \left(1
- 3 \theta_3(\tau)^4  - E_2(\tau)
+ 24 \sum_{n=1}^\infty \Exp{\pi\I  (n^2+n) \tau} \frac{1+\Exp{2\pi\I n\tau}} {(1-\Exp{2\pi\I n\tau})^2}\right).\\
\end{align*}
\end{prop}
\begin{proof}
We use the following decomposition of $\mu$:
\begin{multline*}
\mu(s, z;\tau) = \frac{\Exp{\pi\I s}}{\theta(z; \tau)(1-\Exp{2\pi\I s})} \\
+ \frac{1}{\theta(z;\tau)} \sum_{n=1}^\infty (-1)^n \Exp{\pi\I (n^2+n) \tau} \left( \frac{\Exp{2\pi\I n z + \pi\I s}}{1-\Exp{2\pi\I n\tau+2\pi\I s}} - \frac{\Exp{-2\pi\I n z -\pi\I s}}{1-\Exp{2\pi\I n\tau-2\pi\I s}} \right).
\end{multline*}
 We compute the Taylor expansion with respect to $2\pi \I s$ around
 $s=0$ of the expression above for the following values of $z$:
 $\frac{1+\tau}2$, $\frac{\tau}2$, $\frac12$. We need only the
 coefficient at $2\pi\I s$. This coefficient equals
\[
\frac{1}{24\,\theta(z;\tau)} + \frac{1}{\theta(z;\tau)} \sum_{n=1}^\infty (-1)^n \Exp{\pi\I(n^2+n)\tau} \frac{(1+\Exp{2\pi\I n\tau}) (\Exp{2\pi\I n z} + \Exp{-2\pi\I n z})}{2(1-\Exp{2\pi\I n\tau})^2}.
\]
Therefore in the case $z = \frac{1+\tau}2$ we obtain 
\begin{multline*}
-\frac{\Exp{\frac{\pi\I\tau}4}}{\theta_1(\tau)} \left( 
\frac{1}{24}  + \sum_{n=1}^\infty \Exp{\pi\I n^2 \tau} \frac{(1+\Exp{2\pi\I n\tau}) ^2}{2(1-\Exp{2\pi\I n\tau})^2} \right)
\\=
-\frac{\Exp{\frac{\pi\I\tau}4}}{\theta_1(\tau)} \left(
-\frac{5}{24} + \frac{\theta_1(\tau)}{4} + 
2\sum_{n=1}^\infty \frac{\Exp{\pi\I (n^2+2n)\tau}}{(1-\Exp{2\pi\I n\tau})^2} \right),
\end{multline*}
similarly, in the case $z = \frac{\tau}2$
\[
\I \frac{\Exp{\frac{\pi\I\tau}4}}{\theta_2(\tau)} \left(
-\frac{5}{24} + \frac{\theta_2(\tau)}{4} + 
2\sum_{n=1}^\infty \frac{(-1)^n \Exp{\pi\I (n^2+2n)\tau}}{(1-\Exp{2\pi\I n\tau})^2} \right),
\]
and finally, in the case $z = \frac12$:
\[
-\frac{1}{\theta_3(\tau)} \left(\frac{1}{24} + \sum_{n=1}^\infty \Exp{\pi\I  (n^2+n) \tau} \frac{1+\Exp{2\pi\I n\tau}} {(1-\Exp{2\pi\I n\tau})^2} \right).
\]

For applying Proposition \ref{prop:my3} we also need to compute the coefficients at $2\pi\I s$ of the following expressions (we omit $\tau$ from the arguments of $\theta$):
\[
\frac{\I\eta(\tau)^3 \theta(\tfrac12+s) \theta(\tfrac{\tau}2)}
{\theta(\tfrac{1-\tau}2) \theta(\tfrac12) \theta(s) \theta(\tfrac{\tau}2+s)},
\qquad
\frac{\I\eta(\tau)^3 \theta(\tfrac12+s) \theta(\tfrac{1+\tau}2)}
{\theta(-\tfrac{\tau}2) \theta(\tfrac12) \theta(s) \theta(\tfrac{1+\tau}2+s)},
\]
\[
\frac{\I\eta(\tau)^3 \theta(\tfrac{\tau}2+s) \theta(\tfrac{1+\tau}2)}
{\theta(-\tfrac12) \theta(\tfrac{\tau}2) \theta(s) \theta(\tfrac{1+\tau}2+s)},
\]
For this we need to compute the Taylor expansions of $\theta(s)$, $\theta(\tfrac12+s)$, $\theta(\tfrac{\tau}2+s)$ and $\theta(\tfrac{1+\tau}2+s)$ up to second term with respect to $2\pi\I s$. We have (denoting by $'$ the operator $\frac{d}{2\pi\I d\tau}$):
\begin{align*}
\theta(s) 
&= (2\pi\I s) \I \eta(\tau)^3 \left(1 + (2\pi\I s)^2 \frac{\eta'(\tau)}{\eta(\tau)}\right) + \cdots,\\
\theta(\tfrac12+s) 
&= -\theta_3(\tau) \left(1 + (2\pi\I s)^2 \frac{\theta_3'(\tau)}{\theta_3(\tau)}\right) + \cdots,\\
\theta(\tfrac{\tau}2+s)
&= -\I \Exp{-\frac{\pi\I\tau}4} \theta_2(\tau) \left(1 - \frac{2\pi\I s}{2} + (2\pi\I s)^2\left(\frac18+\frac{\theta_2'(\tau)}{\theta_2(\tau)}\right)\right) + \cdots,\\
\theta(\tfrac{1+\tau}2+s)
&= -\Exp{-\frac{\pi\I\tau}4} \theta_1(\tau) \left(1 - \frac{2\pi\I s}{2} + (2\pi\I s)^2\left(\frac18+\frac{\theta_1'(\tau)}{\theta_1(\tau)}\right)\right) + \cdots.
\end{align*}

Thus the coefficients at $2\pi\I s$ of the expressions in question are, correspondingly,
\begin{align*}
-\frac{\Exp{\frac{\pi\I\tau}4}}{\theta_1(\tau)} & \left(
\frac{\theta_3'(\tau)}{\theta_3(\tau)} - \frac{\eta'(\tau)}{\eta(\tau)} + \frac18 - \frac{\theta_2'(\tau)}{\theta_2(\tau)} \right), \\
-\I\frac{\Exp{\frac{\pi\I\tau}4}}{\theta_2(\tau)} & \left(
\frac{\theta_3'(\tau)}{\theta_3(\tau)} - \frac{\eta'(\tau)}{\eta(\tau)} + \frac18 - \frac{\theta_1'(\tau)}{\theta_1(\tau)} \right), \\
\frac{1}{\theta_3(\tau)} & \left(
\frac{\theta_2'(\tau)}{\theta_2(\tau)} - \frac{\eta'(\tau)}{\eta(\tau)} - \frac{\theta_1'(\tau)}{\theta_1(\tau)} \right).
\end{align*}

Putting everything together 
\begin{align*}
&\frac{d}{2\pi\I d s}\bigg\vert_{s=0} \mu\left(\frac12, \frac{\tau}2+s;\tau\right)\\
&= \frac{d}{2\pi\I d s}\bigg\vert_{s=0} \left( \mu\left(\frac{1-\tau} 2, s;\tau\right) + \frac{\I\eta(\tau)^3 \theta(\tfrac12+s) \theta(\tfrac{\tau}2)}
{\theta(\tfrac{1-\tau}2) \theta(\tfrac12) \theta(s) \theta(\tfrac{\tau}2+s)} \right)\\
&=-\frac{\Exp{\frac{\pi\I\tau}4}}{\theta_1(\tau)} \left(
-\frac{1}{12} + \frac{\theta_1(\tau)}{4} + 
\frac{\theta_3'(\tau)}{\theta_3(\tau)} - \frac{\eta'(\tau)}{\eta(\tau)}  - \frac{\theta_2'(\tau)}{\theta_2(\tau)} + 2\sum_{n=1}^\infty \frac{\Exp{\pi\I (n^2+2n)\tau}}{(1-\Exp{2\pi\I n\tau})^2}\right),
\end{align*}
\begin{align*}
&\frac{d}{2\pi\I d s}\bigg\vert_{s=0} \mu\left(\frac12, \frac{1+\tau}2+s;\tau\right)\\
&= \frac{d}{2\pi\I d s}\bigg\vert_{s=0} \left( \mu\left(-\frac{\tau} 2, s;\tau\right) + 
\frac{\I\eta(\tau)^3 \theta(\tfrac12+s) \theta(\tfrac{1+\tau}2)}
{\theta(-\tfrac{\tau}2) \theta(\tfrac12) \theta(s) \theta(\tfrac{1+\tau}2+s)}\right)\\
&=-\I \frac{\Exp{\frac{\pi\I\tau}4}}{\theta_2(\tau)} \left(
-\frac{1}{12} + \frac{\theta_2(\tau)}{4} + 
\frac{\theta_3'(\tau)}{\theta_3(\tau)} - \frac{\eta'(\tau)}{\eta(\tau)} - \frac{\theta_1'(\tau)}{\theta_1(\tau)}
+2\sum_{n=1}^\infty \frac{(-1)^n \Exp{\pi\I (n^2+2n)\tau}}{(1-\Exp{2\pi\I n\tau})^2} \right),
\end{align*}
\begin{align*}
&\frac{d}{2\pi\I d s}\bigg\vert_{s=0} \mu\left(\frac{\tau}2, \frac{1+\tau}2+s;\tau\right)\\
&= \frac{d}{2\pi\I d s}\bigg\vert_{s=0} \left( \mu\left(-\frac12, s;\tau\right) + 
\frac{\I\eta(\tau)^3 \theta(\tfrac{\tau}2+s) \theta(\tfrac{1+\tau}2)}
{\theta(-\tfrac12) \theta(\tfrac{\tau}2) \theta(s) \theta(\tfrac{1+\tau}2+s)}\right)\\
&=\frac{1}{\theta_3(\tau)} \left(\frac{1}{24} 
+ \frac{\theta_2'(\tau)}{\theta_2(\tau)} - \frac{\eta'(\tau)}{\eta(\tau)} - \frac{\theta_1'(\tau)}{\theta_1(\tau)} 
+ \sum_{n=1}^\infty \Exp{\pi\I  (n^2+n) \tau} \frac{1+\Exp{2\pi\I n\tau}} {(1-\Exp{2\pi\I n\tau})^2}\right).
\end{align*}
The statements we need to prove follow from the following identities:
\begin{align*}
\frac{\theta_3'(\tau)}{\theta_3(\tau)} - \frac{\eta'(\tau)}{\eta(\tau)}  - \frac{\theta_2'(\tau)}{\theta_2(\tau)}
&= \frac{\theta_1^4(\tau)}8 - \frac{E_2(\tau)}{24},\\
\frac{\theta_3'(\tau)}{\theta_3(\tau)} - \frac{\eta'(\tau)}{\eta(\tau)}  - \frac{\theta_1'(\tau)}{\theta_1(\tau)}
&= \frac{\theta_2^4(\tau)}8 - \frac{E_2(\tau)}{24},\\
\frac{\theta_2'(\tau)}{\theta_2(\tau)} - \frac{\eta'(\tau)}{\eta(\tau)}  - \frac{\theta_1'(\tau)}{\theta_1(\tau)}
&= -\frac{\theta_3^4(\tau)}8 - \frac{E_2(\tau)}{24},\\
\end{align*}
\end{proof}

Propositions \ref{prop:my7} and \ref{prop:my8} together give the Fourier expansions of $\wt h_i$. Denote
\begin{align*}
h_1(\tau) &= \frac{1}{24 \theta_1(\tau)} \left(
2  + E_2(\tau)
 - 48 \sum_{n=1}^\infty \frac{\Exp{\pi\I (n^2+2n)\tau}}{(1-\Exp{2\pi\I n\tau})^2}\right),\\
h_2(\tau) &= \frac{1}{24 \theta_2(\tau)} \left(
2 + E_2(\tau)
 -48 \sum_{n=1}^\infty \frac{(-1)^n \Exp{\pi\I (n^2+2n)\tau}}{(1-\Exp{2\pi\I n\tau})^2} \right),\\
h_3(\tau) &= \frac{1}{24 \theta_3(\tau)} \left(
-1 + E_2(\tau)
 -24 \sum_{n=1}^\infty \Exp{\pi\I  (n^2+n) \tau} \frac{1+\Exp{2\pi\I n\tau}} {(1-\Exp{2\pi\I n\tau})^2}\right).
\end{align*}
The series $h_i$ are holomorphic power series converging on the upper half plane. Denote
\begin{align*}
R_1(\tau) &= \frac12 \sum_{n\in\Z} |n| \beta(2y n^2)  \Exp{ - \pi\I n^2\tau}
 - \frac{\ol{\theta_1(\tau)}}{2\pi\sqrt{2y}},\\
R_2(\tau) &= \frac12 \sum_{n\in\Z} (-1)^n |n| \beta(2y n^2)  \Exp{ - \pi\I n^2\tau}
 - \frac{\ol{\theta_2(\tau)}}{2\pi\sqrt{2y}},\\
R_3(\tau) &= \frac12 \sum_{\nu\in\Z+\frac12} |\nu| \beta(2y\nu^2) \Exp{ -\pi\I\nu^2\tau} 
 - \frac{\ol{\theta_3(\tau)}}{2\pi\sqrt{2y}}.
\end{align*}

\begin{prop} 
For $i=1,2,3$ we have
\[
\wt h_i(\tau) = h_i(\tau) + R_i(\tau).
\]
\end{prop}

In his thesis Zwegers also represents $R$ as a certain integral
involving a theta function of weight $\frac32$. In our case, we also
have such a representation but with theta functions of weight $\frac12$.
\begin{prop}\label{prop:my10}
For $i=1,2,3$ 
\[
R_i(\tau) = \frac{1}{4\pi\I} \;\ol{\int_{\tau}^{\I\infty} \frac{\theta_i(z) dz}{(-\I(z-\ol\tau))^{\frac32}}}.
\]
\end{prop}
\begin{proof}
Note that the integral on the right converges. We prove the identity termwise using the following formula for a real number $a$:
\begin{align*}
\frac{1}{2} |a| \beta(2 y a^2) \Exp{-\pi\I a^2\tau} - \frac{\ol{\Exp{\pi\I a^2\tau}}}{2\pi\sqrt{2y}}
&= \ol{-\int_\tau^{\I\infty} \frac{a^2 \I \Exp{\pi\I a^2 z} dz}{2 \sqrt{-\I(z-\ol\tau)}} - \frac{\Exp{\pi\I a^2\tau}}{2\pi\sqrt{2y}} } \\
&= \frac{1}{4\pi\I}\; \ol{\int_{\tau}^{\I\infty} \frac{\Exp{\pi\I a^2 z} dz} {(-\I(z-\ol\tau))^{\frac32}} }.
\end{align*}
This formula is obtained using the integral representation of $\beta$
\[
\beta(2 y a^2) = \int_{2 y a^2}^\infty \Exp{-\pi t}\frac{dt}{\sqrt{t}}
\]
after the substitution $t = -\I(z-\ol\tau) a^2$, and then integration by parts. The case $a=0$ should be considered separately.
\end{proof}

We would like to compute the Fourier coefficients of $h_i\theta_i$
explicitly. Since we know the Fourier coefficients of $E_2$, it remains
to consider the following expressions:
\begin{prop}
We have
\begin{align*}
2 \sum_{n=1}^\infty \frac{\Exp{\pi\I (n^2+2n)\tau}}{(1-\Exp{2\pi\I n\tau})^2}
&= \sum_{\substack{m>n>0, \\m-n\;\mathrm{even}}} m \Exp{\pi\I m n \tau} - 
\sum_{\substack{n>m>0, \\m-n\;\mathrm{even}}} m \Exp{\pi\I m n \tau} \\
2 \sum_{n=1}^\infty \frac{(-1)^n \Exp{\pi\I (n^2+2n)\tau}}{(1-\Exp{2\pi\I n\tau})^2}
&= \sum_{\substack{m>n>0, \\m-n\;\mathrm{even}}} m (-1)^m \Exp{\pi\I m n \tau} - 
\sum_{\substack{n>m>0, \\m-n\;\mathrm{even}}} m (-1)^m \Exp{\pi\I m n \tau} \\
\sum_{n=1}^\infty \Exp{\pi\I  (n^2+n) \tau} \frac{1+\Exp{2\pi\I n\tau}} {(1-\Exp{2\pi\I n\tau})^2}
&= \sum_{\substack{m>n>0, \\m-n\;\mathrm{odd}}} m \Exp{\pi\I m n \tau} - 
\sum_{\substack{n>m>0, \\m-n\;\mathrm{odd}}} m \Exp{\pi\I m n \tau} \\
\end{align*}
\end{prop}
\begin{proof}
It is clear.
\end{proof}

Looking at the expansions we observe that
\begin{multline}
2 \sum_{n=1}^\infty \frac{\Exp{\pi\I (n^2+2n)\tau}}{(1-\Exp{2\pi\I n\tau})^2}
+ 2 \sum_{n=1}^\infty \frac{(-1)^n \Exp{\pi\I (n^2+2n)\tau}}{(1-\Exp{2\pi\I n\tau})^2}
\\= 4\left(\sum_{m>n>0} m \Exp{4 \pi\I m n \tau} - 
\sum_{n>m>0} m \Exp{4 \pi\I m n \tau}\right),
\end{multline}
and
\begin{multline}\label{eq:my30}
2 \sum_{n=1}^\infty \frac{\Exp{\pi\I (n^2+2n)\tau}}{(1-\Exp{2\pi\I n\tau})^2}
+ \sum_{n=1}^\infty \Exp{\pi\I  (n^2+n) \tau} \frac{1+\Exp{2\pi\I n\tau}} {(1-\Exp{2\pi\I n\tau})^2}
\\= \sum_{m>n>0} m \Exp{\pi\I m n \tau} - 
\sum_{n>m>0} m \Exp{\pi\I m n \tau}
\end{multline}

Therefore it is not difficult to complete the proof of the following statement:
\begin{prop}
We have
\[
h_1(\tau) \theta_1(\tau) + h_2(\tau) \theta_2(\tau) - 4 \Bigl(h_1(4 \tau) \theta_1(4 \tau) + h_3(4 \tau) \theta_3(4 \tau)\Bigr) = -\frac{\theta_1(2\tau)^4}4.
\]
\end{prop}

Using the corresponding identities between the theta functions, namely
\[
\theta_1(\tau)+\theta_2(\tau) = 2 \theta_1(4\tau),\qquad
\theta_1(\tau)+\theta_3(\tau) = \theta_1(\tfrac\tau4),
\]
the integral representation of $R_i$ from Proposition \ref{prop:my10} and the change of variables
\[
\ol{4\pi\I R_i(4\tau)} = \int_{4 \tau}^{\I\infty} \frac{\theta_i(z) dz}{(-\I(z-4\ol\tau))^{\frac32}}
= \frac12 \int_{\tau}^{\I\infty} \frac{\theta_i(4 z) dz}{(-\I(z-\ol\tau))^{\frac32}}
\]
we obtain
\[
R_1(\tau) \theta_1(\tau) + R_2(\tau) \theta_2(\tau) - 4 \left(R_1(4\tau) \theta_1(4\tau) + R_3(4\tau) \theta_3(4\tau)\right) = 0.
\]
Therefore we also have the following.
\begin{prop}
\[
\wt h_1(\tau) \theta_1(\tau) + \wt h_2(\tau) \theta_2(\tau) - 4 \left(\wt h_1(4 \tau) \theta_1(4 \tau) + \wt h_3(4 \tau) \theta_3(4 \tau)\right) = -\frac{\theta_1(2\tau)^4}4.
\]
\end{prop}

Now we are ready to prove Theorem 1.5. This is done in a series of propositions. Take
\[
\wt g(\tau) = -\frac{\wt h_1(2\tau)\theta_1(2\tau) + \wt h_3(2\tau)\theta_3(2\tau)}2 + \frac{\theta_1(\tau)^4 + \theta_2(\tau)^4}{96},
\]
then the following holds.
\begin{prop}\label{prop:my15}
The function $\wt g(\tau)$ transforms like a modular form of weight $2$:
\[
\wt g(\tau+1) = \wt g(\tau), \qquad \wt g\left(-\frac1\tau\right) = \tau^2 \wt g(\tau).
\]
\end{prop}
The function $\wt g(\tau)$ decomposes as
\[
\wt g(\tau) = g(\tau) + r(\tau),
\]
where
\[
g(\tau) = -\frac{h_1(2\tau)\theta_1(2\tau) + h_3(2\tau)\theta_3(2\tau)}2 + \frac{\theta_1(\tau)^4 + \theta_2(\tau)^4}{96}
\]
and 
\[
r(\tau) = -\frac{R_1(2\tau)\theta_1(2\tau) + R_3(2\tau)\theta_3(2\tau)}2.
\]
It is not difficult to compute the Fourier expansion of $g$:
\begin{prop}\label{prop:my20}
We have
\[
g(\tau) = -\frac{E_2(\tau)}{24} -\frac12 \sum_{n\in \Z\setminus \{0\}} \frac{n q^{n^2}}{1-q^n}
= -\frac{1}{24} + \sum_{n=1}^\infty \sigma'(n) q^n,
\]
where $\sigma'(n)$ denotes the sum of positive divisors of $n$ which are greater than $\sqrt{n}$, plus half $\sqrt{n}$  in the case if $n$ is a perfect square.
\end{prop}
\begin{proof}
Using the third identity from \eqref{eq:my8} and \eqref{eq:my30} we find
\begin{align*}
g(\tau) &= \frac{1}{48}\left(-1 - E_2(\tau) + 24\left(\sum_{m>n>0} m \Exp{2\pi\I m n \tau} - 
\sum_{n>m>0} m \Exp{2\pi\I m n \tau}\right)\right)\\
&=-\frac{1}{24} +\frac12 \sum_{m>0,\; n>0} m \Exp{2\pi\I m n \tau} +\frac12 \sum_{m>n>0} m \Exp{2\pi\I m n \tau} -
\frac12 \sum_{n>m>0} m \Exp{2\pi\I m n \tau}\\
&=-\frac{1}{24} + \sum_{m>0,\; n>0} m \Exp{2\pi\I m n \tau} - \sum_{n\geq m>0} m \Exp{2\pi\I m n \tau} + \frac12 \sum_{n>0} n \Exp{2\pi\I n^2\tau}.
\end{align*}
Then the statement follows.
\end{proof}
Now, by Proposition \ref{prop:my10}, $r(\tau) = R(\overline{D}(r(\tau)))
= R(\overline{D}(\tilde{g}(\tau)))$, where $R$ is the operator from the
introduction, hence $g$ is the holomorphic part of
$\tilde{g}$. Proposition \ref{prop:my20} gives the Fourier expansion of
$g$, which, as one can easily verify, coincides with $\cJ_{-2}$. Having
the transformation properties of $\tilde{g}$ proved in Proposition
\ref{prop:my15}, it remains to check only the uniqueness statement. This
is obvious since there are no holomorphic modular forms of weight $2$
for $\SL(2,\Z)$.

\section{Acknowledgments}
Authors thank Max-Planck-Institut f\"ur Mathematik and Institut 
des Hautes \'Etudes Scientifiques.  They thank Aspinwall, Bridgeland,
Harder, Kajiura, Kontsevich, Mizuno, Takahashi, Toda, Yoshioka, 
Zagier, and Zudilin for their helpful comments or discussions. 
Especially, Kontsevich and Zagier gave them their preprints 
\cite{KoSo} and \cite{Za_Bu} and made useful discussions and remarks.  
The second author thanks University of Texas Austin for its 
support and active sessions with Ben-Zvi, Freed, Keel, and 
Soibelman during his visit in November 07.

\end{document}